\documentclass[11pt,reqno]{amsart}
\usepackage{amssymb,latexsym}
\usepackage{graphicx}
\usepackage{float}
\usepackage{url}		

\calclayout

\makeatletter
\newcommand{\monthyear}[1]{%
  \def\@monthyear{\uppercase{#1}}}
\newcommand{\volnumber}[1]{%
  \def\@volnumber{\uppercase{#1}}}
\AtBeginDocument{%
\def\ps@plain{\ps@empty
  \def\@oddfoot{\@monthyear \hfil \thepage}%
  \def\@evenfoot{\thepage \hfil \@volnumber}}
\def\ps@firstpage{\ps@plain}
\def\ps@headings{\ps@empty
  \def\@evenhead{%
    \setTrue{runhead}%
    \def\thanks{\protect\thanks@warning}%
    \uppercase{The Fibonacci Quarterly}\hfil}%
  \def\@oddhead{%
    \setTrue{runhead}%
    \def\thanks{\protect\thanks@warning}%
    \hfill\uppercase{Hypergeometric Template}}%
  \let\@mkboth\markboth
  \def\@evenfoot{%
    \thepage \hfil \@volnumber}%
  \def\@oddfoot{%
    \@monthyear \hfil \thepage}%
  }%
\footskip=25pt
\pagestyle{headings}%
}
\makeatother

\theoremstyle{plain}
\numberwithin{equation}{section}
\newtheorem{thm}{Theorem}[section]
\newtheorem{theorem}[thm]{Theorem}
\newtheorem{lemma}[thm]{Lemma}
\newtheorem{example}[thm]{Example}

\begin{document}
\monthyear{Month Year}
\volnumber{Volume, Number}
\setcounter{page}{1}

\title{On the $x$--coordinates of Pell equations which are products of two Lucas numbers}
\author{Mahadi Ddamulira}
\address{Institute of Analysis and Number Theory\\ 
        Graz University of Technology\\ \newline \indent
		Kopernikusgasse 41/II\\ 
		A-8010 Graz, Austria}
\email{mddamulira@tugraz.at; mahadi@aims.edu.gh}
\thanks{This research was supported by the Austrian Science Fund (FWF) grants: F5510-N26 -- Part of the special research program (SFB), ``Quasi-Monte Carlo Methods: Theory and Applications'' and W1230 --``Doctoral Program Discrete Mathematics''.}

\begin{abstract}
Let $ \{L_n\}_{n\ge 0} $ be the sequence of Lucas numbers given by $ L_0=2, ~ L_1=1 $ and $ L_{n+2}=L_{n+1}+L_n $ for all $ n\ge 0 $. In this paper, for an integer $d\geq 2$ which is square-free, we show that there is at most one value of the positive integer $x$ participating in the Pell equation $x^{2}-dy^{2}
=\pm 1$ which is a product of two Lucas numbers, with a few exceptions that we completely characterize.
\end{abstract} 

\maketitle

\section{Introduction}
Let $ \{L_n\}_{n\ge 0} $ be the sequence of Lucas numbers given by $ L_0=2, ~ L_1=1 $ and $$ L_{n+2}=L_{n+1}+L_n $$ for all $ n\ge 0 $. This is sequence A000032 on the Online Encyclopedia of Integer Sequences (OEIS). The first few terms of this sequence are
$$ \{L_n\}_{n\ge 0} \quad =\quad 2, 1, 3, 4, 7, 11, 18, 29, 47, 76, 123, 199, 322, 521, 843, 1364, 2207, 3571, \ldots. $$
Putting $ \displaystyle{(\alpha, \beta) = \left(\frac{1+\sqrt{5}}{2}, \frac{1-\sqrt{5}}{2}\right)} $ for the roots of the characteristic equation  $ \displaystyle{r^2-r-1=0} $ of the Lucas sequence, the Binet formula for its general terms is given by
\begin{eqnarray}\label{Binet}
L_n = \alpha^{n}+\beta^{n}, \qquad \text{ for all } \quad n\ge 0.
\end{eqnarray}
Furthermore, we can prove by induction that the inequality
\begin{eqnarray}\label{Pell11}
\alpha^{n-1}\le L_n\le \alpha^{n+2},
\end{eqnarray}
holds for all $ n\ge 0 $.

Let $d\geq 2$ be a positive integer which is not a perfect square. It is well known that the Pell equation
\begin{eqnarray}
x^{2}-dy^{2}=\pm 1 \label{Pelleqn}
\end{eqnarray}
has infinitely many positive integer solutions $(x,y)$. By putting $(x_1, y_1)$ for the smallest positive solution, all solutions are of the form $ (x_k, y_k) $ for some positive integer $k$, where
\begin{eqnarray}
x_k+y_k\sqrt{d} = (x_1+y_1\sqrt{d})^n\qquad {\text{\rm for~all}} \quad k\ge 1.\label{Pellsoln}
\end{eqnarray}
Furthermore, the sequence $ \{x_k\}_{k\ge 1} $ is binary recurrent. In fact, the following formula
\begin{eqnarray*}
x_k=\dfrac{(x_1+y_1\sqrt{d})^{k}+(x_1-y_1\sqrt{d})^{k}}{2},
\end{eqnarray*}
holds for all positive integers $ k $.

Recently, Kafle et al. \cite{Bir1} considered the Diophantine equation
\begin{equation}\label{eq:Luca1} 
x_{n} =F_{\ell}F_{m},
\end{equation}
where $\{F_m\}_{m \ge 0}$ is the sequence of Fibonacci numbers given
by $F_0 = 0$, $F_1 = 1$ and $F_{m+2} = F_{m+1} + F_m$ for all $m \ge 0$. They proved that  equation \eqref{eq:Luca1} has at most one solution $n$ in positive integers except for $d=2, 3, 5$, for which case equation \eqref{eq:Luca1} has the  solutions $ x_1=1$ and $ x_2=3 $, $ x_1=2 $ and $ x_2=26 $, $ x_1=2$ and $x_2=9 $, respectively.

There are many other researchers who have studied related problems involving the intersection sequence $ \{x_n\}_{n\ge 1} $ with linear recurrence sequences of interest. For example, see \cite{Bravo, Ddamulira, Ddamulira1, Ddamulira2, BLT1, BLT2, BLT3, Luca16, Luca15, Togbe}.

\section{Main Result}

In this paper, we study a similar problem to that of Kafle et al. \cite{Bir1}, but with the Lucas numbers instead of the Fibonacci numbers. That is, we show that there is at most one value of the positive integer $x$ participating in \eqref{Pelleqn} which is a product of two Lucas numbers, with a few exceptions that we completely cahracterize. This can be interpreted as solving the Diophantine  equation
\begin{align}\label{Problem}
 x_{k}=L_{n}L_{m},
\end{align}
in nonnegative integers $ (k, n, m) $ with $ k\ge 1 $ and $ 0\le m \le  n $.

\begin{theorem}\label{Main}
For each square-free integer $ d\ge 2 $ there is at most one integer $ k $ such that the equation \eqref{Problem} holds, except for $ d\in \{2, 3, 5,15,17, 35\} $ for which $ x_1=1, ~ x_2=3, x_3=7, x_9= 1393  $ (for $ d=2 $), $ x_1=2, ~x_2=7 $ (for $ d=3 $), $ x_1=2, ~ x_2=9 $ (for $ d=5 $), $ x_1=4, ~ x_5=15124 $ (for $ d=15 $), $ x_1=4, ~x_2=33 $ (for $ d=17 $) and $ x_1=6, ~x_3=846 $ (for $ d=35 $).
\end{theorem}

\section{Preliminary Results}

\subsection{Notations and terminology from algebraic number theory} 

We begin by recalling some basic notions from algebraic number theory.

Let $\eta$ be an algebraic number of degree $d$ with minimal primitive polynomial over the integers
$$
a_0x^{d}+ a_1x^{d-1}+\cdots+a_d = a_0\prod_{i=1}^{d}(x-\eta^{(i)}),
$$
where the leading coefficient $a_0$ is positive and the $\eta^{(i)}$'s are the conjugates of $\eta$. Then the \textit{logarithmic height} of $\eta$ is given by
$$ 
h(\eta):=\dfrac{1}{d}\left( \log a_0 + \sum_{i=1}^{d}\log\left(\max\{|\eta^{(i)}|, 1\}\right)\right).
$$
In particular, if $\eta=p/q$ is a rational number with $\gcd (p,q)=1$ and $q>0$, then $h(\eta)=\log\max\{|p|, q\}$. The following are some of the properties of the logarithmic height function $h(\cdot)$, which will be used in the next sections of this paper without reference:
\begin{eqnarray}
h(\eta\pm \gamma) &\leq& h(\eta) +h(\gamma) +\log 2,\nonumber\\
h(\eta\gamma^{\pm 1})&\leq & h(\eta) + h(\gamma),\\
h(\eta^{s}) &=& |s|h(\eta) \qquad (s\in\mathbb{Z}). \nonumber
\end{eqnarray}

\subsection{Linear forms in logarithms}
In order to prove our main result Theorem \ref{Main}, we need to use several times a Baker--type lower bound for a nonzero linear form in logarithms of algebraic numbers. There are many such 
in the literature  like that of Baker and W{\"u}stholz from \cite{bawu07}.  We start by recalling the result of Bugeaud, Mignotte and Siksek (\cite{BuMiSi}, Theorem 9.4, pp. 989), which is a  modified version of the result of Matveev \cite{MatveevII}, which is one of our main tools in this paper.

\begin{theorem}\label{Matveev11} Let $\gamma_1,\ldots,\gamma_t$ be positive real  numbers in a number field 
$\mathbb{K} \subseteq  \mathbb{R}$ of degree $D$, $b_1,\ldots,b_t$ be nonzero integers, and assume that
\begin{equation}
\label{eq:Lambda}
\Lambda:=\gamma_1^{b_1}\cdots\gamma_t^{b_t} - 1,
\end{equation}
is nonzero. Then
$$
\log |\Lambda| > -1.4\times 30^{t+3}\times t^{4.5}\times D^{2}(1+\log D)(1+\log B)A_1\cdots A_t,
$$
where
$$
B\geq\max\{|b_1|, \ldots, |b_t|\},
$$
and
$$A
_i \geq \max\{Dh(\gamma_i), |\log\gamma_i|, 0.16\},\qquad {\text{for all}}\qquad i=1,\ldots,t.
$$
\end{theorem}
When $t=2$  and  $\gamma_1, ~\gamma_2$ are positive and multiplicatively independent, we can use a result of Laurent, Mignotte and Nesterenko \cite{Laurent:1995}. Namely, let in this case $B_1, ~B_2$ be real numbers larger than $1$ such that
\begin{eqnarray*}
\log B_i\geq \max\left\{h(\gamma_i), \dfrac{|\log\gamma_i|}{D}, \dfrac{1}{D}\right\},\qquad {\text{\rm for}}\quad  i=1,2,
\end{eqnarray*}
and put
\begin{eqnarray*}
b^{\prime}:=\dfrac{|b_{1}|}{D\log B_2}+\dfrac{|b_2|}{D\log B_1}.
\end{eqnarray*}
Put
\begin{equation}
\label{eq:Gamma}
\Gamma:= b_1\log\gamma_1+b_2\log\gamma_2.
\end{equation}
We note that $\Gamma\neq0$ because $\gamma_{1} $ and $\gamma_{2}$ are multiplicatively independent. The following result is Corollary $ 2 $ in \cite{Laurent:1995}.
\begin{theorem}\label{Matveev12}
With the above notations, assuming that $ \eta_{1}, \eta_{2}$ are positive and multiplicatively independent, then
\begin{eqnarray}
\log |\Gamma|> -24.34D^4\left(\max\left\{\log b^{\prime}+0.14, \dfrac{21}{D}, \dfrac{1}{2}\right\}\right)^{2}\log B_1\log B_2.
\end{eqnarray}
\end{theorem}
Note that with $\Gamma$ given by \eqref{eq:Gamma}, we have $e^{\Gamma}-1=\Lambda$, where $\Lambda$ is given by \eqref{eq:Lambda} in case $t=2$, which explains the connection between Theorem \ref{Matveev11} and Theorem \ref{Matveev12}.

\subsection{Reduction procedure}\label{Reduction}
During the calculations, we get upper bounds on our variables which are too large, thus we need to reduce them. To do so, we use some results from the theory of continued fractions. 

For the treatment of linear forms homogeneous in two integer variables, we use the well-known classical result in the theory of Diophantine approximation.
\begin{lemma}\label{Legendre}
Let $\tau$ be an irrational number,  $ \frac{p_0}{q_0}, \frac{p_1}{q_1}, \frac{p_2}{q_2}, \ldots $ be all the convergents of the continued fraction of $ \tau $ and $ M $ be a positive integer. Let $ N $ be a nonnegative integer such that $ q_N> M $. Then putting $ a(M):=\max\{a_{i}: i=0, 1, 2, \ldots, N\} $, the inequality
\begin{eqnarray*}
\left|\tau - \dfrac{r}{s}\right|> \dfrac{1}{(a(M)+2)s^{2}},
\end{eqnarray*}
holds for all pairs $ (r,s) $ of positive integers with $ 0<s<M $.
\end{lemma}

For a nonhomogeneous linear form in two integer variables, we use a slight variation of a result due to Dujella and Peth{\H o} (see \cite{dujella98}, Lemma 5a). For a real number $X$, we write  $||X||:= \min\{|X-n|: n\in\mathbb{Z}\}$ for the distance from $X$ to the nearest integer.
\begin{lemma}\label{Dujjella}
Let $M$ be a positive integer, $\frac{p}{q}$ be a convergent of the continued fraction of the irrational number $\tau$ such that $q>6M$, and  $A,B,\mu$ be some real numbers with $A>0$ and $B>1$. Let further 
$\varepsilon: = ||\mu q||-M||\tau q||$. If $ \varepsilon > 0 $, then there is no solution to the inequality
$$
0<|u\tau-v+\mu|<AB^{-w},
$$
in positive integers $u,v$ and $w$ with
$$ 
u\le M \quad {\text{and}}\quad w\ge \dfrac{\log(Aq/\varepsilon)}{\log B}.
$$
\end{lemma}

At various occasions, we need to find a lower bound for linear forms in logarithms with bounded integer coefficients in three and four variables. In this case we use the LLL  algorithm that we describe below. Let $ \tau_1, \tau_2, \ldots \tau_t \in\mathbb{R}$ and the linear form
\begin{eqnarray}
x_1\tau_1+x_2\tau_2+\cdots+x_t\tau_t \quad \text{ with } \quad |x_i|\le X_i.
\end{eqnarray}
We put $ X:=\max\{X_i\} $, $ C> (tX)^{t} $ and consider the integer lattice $ \Omega $ generated by
\begin{eqnarray*}
\textbf{b}_j: = \textbf{e}_j+\lfloor C\tau_j\rceil \quad \text{ for} \quad 1\le j\le t-1 \quad \text{ and} \quad \textbf{b}_t:=\lfloor C\tau_t\rceil \textbf{e}_t,
\end{eqnarray*}
where $ C $ is a sufficiently large positive constant.
\begin{lemma}\label{LLL}
Let $ X_1, X_2, \ldots, X_t $ be positive integers such that $ X:=\max\{X_i\} $ and $ C> (tX)^{t} $ is a fixed sufficiently large constant. With the above notation on the lattice $ \Omega $, we consider a reduced base $ \{\textbf{b}_i \}$ to $ \Omega $ and its associated Gram-Schmidt orthogonalization base $ \{\textbf{b}_i^*\}$. We set
\begin{eqnarray*}
c_1:=\max_{1\le i\le t}\dfrac{||\textbf{b}_1||}{||\textbf{b}_i^*||}, \quad \theta:=\dfrac{||\textbf{b}_1||}{c_1}, \quad Q:=\sum_{i=1}^{t-1}X_i^{2} \quad \text{and} \quad R:=\left(1+ \sum_{i=1}^{t}X_i\right)/2.
\end{eqnarray*}
If the integers $ x_i $ are such that $ |x_i|\le X_i $, for $ 1\le i \le t $ and $ \theta^2\ge Q+R^2 $, then we have
\begin{eqnarray*}
\left|\sum_{i=1}^{t}x_i\tau_i\right|\ge \dfrac{\sqrt{\theta^2-Q}-R}{C}.
\end{eqnarray*}
\end{lemma}
\noindent
For the proof and further details, we refer the reader to the book of Cohen. (Proposition 2.3.20 in \cite{Cohen}, pp. 58--63).

\subsection{Pell equations and Dickson polynomials}
\label{subs:Pell}
Here we give some relations about Pell equations and Dickson polynomials that will be useful in the next section of this paper.

Let $d\ge 2$ be a squarefree integer.  We put $\delta:=x_1+{\sqrt{x_1^2-\epsilon}}$ for the smallest positive integer $x_1$ such that
$$
x_1^2-dy_1^2=\epsilon,\qquad \epsilon\in \{\pm 1\}
$$
for some positive integer $y_1$. Then,
$$
x_k+{y_k\sqrt{d}}=\delta^k\qquad {\text{\rm and}}\qquad x_k-{{y_k}\sqrt{d}}=\eta^k,\qquad {\text{\rm where}}\qquad \eta:=\epsilon \delta^{-1}.
$$
From the above, we get
\begin{equation}
\label{eq:Pellsol1}
2x_k=\delta^k+(\epsilon \delta^{-1})^k\qquad {\text{\rm for~all}}\qquad k\ge 1.
\end{equation}
There is a formula expressing $2x_{k}$ in terms of $2x_1$ by means of the Dickson polynomial $D_{k}(2x_1,\epsilon)$, where 
$$
D_k(x,y)=\sum_{i=0}^{\lfloor k/2\rfloor} \frac{k}{k-i} \binom{k-i}{i}(-y)^i x^{k-2i}.
$$
These polynomials appear naturally in many number theory problems and results, for example  in a result of Bilu and Tichy \cite{BiluTichy} concerning 
polynomials $f(X), g(X)\in {\mathbb Z}[X]$ such that the Diophantine equation $f(x)=g(y)$ has infinitely many integer solutions $(x,y)$. 

\begin{example}
\label{exa:2}
\begin{itemize}
\item[(i)] $k=2$. We have
$$
2x_2=\sum_{i=0}^1 \frac{2}{2-i}\binom{2-i}{i} (-\epsilon)^i (2x_1)^{2-2i}=4x_1^2-2\epsilon, \quad {\text{so}}\quad x_2=2x_1^2-\epsilon.
$$
\item[(ii)] $k=3$. We have
$$
2x_3=\sum_{i=0}^1 \frac{3}{3-i}\binom{3-i}{i} (-\epsilon)^i (2x_1)^{3-2i}=(2x_1)^3-3\epsilon (2x_1),\quad {\text{so}}\quad x_3=4x_1^3-3\epsilon x_1.
$$
\end{itemize}
\end{example}

\section{Bounding the variables}
We assume that $ (x_1, y_1) $ is the smallest positive solution of the Pell equation \eqref{Pelleqn}. As in Subsection \ref{subs:Pell}, we set
$$
x_1^2-dy_1^2=:\epsilon, \qquad \epsilon \in\{\pm 1\},
$$
and put
$$
\delta:=x_1+{\sqrt{d}}y_1\qquad {\text{\rm and  }}\qquad \eta:=x_1-{\sqrt{d}}y_1 =\epsilon \delta^{-1}.
$$
From \eqref{Pellsoln}, we get
\begin{equation}\label{Pellxterm}
x_k=\dfrac{1}{2}\left(\delta^{k}+\eta^{k}\right).
\end{equation}
Since $ \delta\geq 1+\sqrt{2}>\alpha^{3/2}$, it follows that the estimate
\begin{eqnarray}\label{Pellenq}
\dfrac{\delta^{k}}{\alpha^2}\leq x_{k}<\dfrac{\delta^k}{\alpha}\quad {\text{\rm holds for all}} \quad k\geq 1.
\end{eqnarray}
We let $ (k, n, m):=(k_i, n_i, m_i) $ for $ i=1,2 $ be the solutions of \eqref{Problem}. By \eqref{Pell11} and \eqref{Pellenq}, we get 
\begin{eqnarray}\label{Pell1n}
\alpha^{n+m-2}\le L_nL_m = x_k < \dfrac{\delta^{k}}{\alpha} \quad \text{ and } \quad  \dfrac{\delta^{k}}{\alpha^2}\le  x_k = L_nL_m \le \alpha^{n+m+4},
\end{eqnarray}
so
\begin{eqnarray}\label{Pell12}
kc_1\log\delta - 6 < n+m < kc_1\log\delta +1 \quad \text{ where } \quad c_1:=\dfrac{1}{\log\alpha}.
\end{eqnarray}
To fix ideas, we assume that 
\begin{eqnarray*}
n\ge m \qquad \text{ and } \qquad k_1< k_2.
\end{eqnarray*}
We also put
\begin{eqnarray*}
m_3:=\min \{m_1, m_2\}, \quad m_4:=\max \{m_1, m_2\}, \quad n_3:=\min \{n_1, n_2\}, \quad n_4:=\max \{n_1, n_2\}.
\end{eqnarray*}
Using the inequality \eqref{Pell12} together with the fact that $ \delta \geq 1+\sqrt{2} = \alpha^{3/2} $ (so, $ c_1\log \delta > 3/2 $), gives us that 
\begin{eqnarray*}
\frac{3}{2}k_2 < k_2c_1\log\delta < 2n_2+6\le 2n_4+6,
\end{eqnarray*}
so
\begin{eqnarray}\label{DMA1}
k_1<k_2< \frac{4}{3}n_4+4.
\end{eqnarray}
Thus, it is enough to find an upper bound on $ n_4 $. Substituting \eqref{Binet} and \eqref{Pellxterm} in \eqref{Problem} we get
\begin{eqnarray}
\dfrac{1}{2}(\delta^{k}+\eta^{k})= (\alpha^{n}+\beta^{n})(\alpha^{m}+\beta^{m}).
\end{eqnarray}
This can be regrouped as
\begin{eqnarray*}
\delta^{k}2^{-1}\alpha^{-n-m}-1 = -2^{-1}\eta^{k}\alpha^{-n-m}+ (\beta \alpha^{-1})^{n}+(\beta \alpha^{-1})^{m}+(\beta \alpha^{-1})^{n+m}.
\end{eqnarray*}
Since $ \beta = -\alpha^{-1} $, $ \eta = \varepsilon \delta^{-1} $ and using the fact that $ \delta^{k}\geq \alpha^{n+m-1} $ (by \eqref{Pell1n}), we get
\begin{eqnarray*}
\left|\delta^{k}2^{-1}\alpha^{-n-m}-1 \right|&\le& \dfrac{1}{2\delta^{k}\alpha^{n+m}}+\dfrac{1}{\alpha^{2n}}+\dfrac{1}{\alpha^{2m}}+\dfrac{1}{\alpha^{2(n+m)}}\\
&\le& \dfrac{\alpha}{2\alpha^{2(n+m)}}+\dfrac{3}{\alpha^{2m}}\quad< \quad \dfrac{6}{\alpha^{2m}},
\end{eqnarray*}
In the above, we have also  used the facts that  $ n\ge m $ and $ (1/2)\alpha +3 < 6$. Hence,
\begin{eqnarray}\label{Pelling11}
\left|\delta^{k}2^{-1}\alpha^{-n-m}-1 \right| < \dfrac{6}{\alpha^{2m}}.
\end{eqnarray} 
We let $\Lambda_1: = \delta^{k}2^{-1}\alpha^{-n-m}-1$. We put 
\begin{equation}
\label{eq:Gamma1}
\Gamma_1:=k\log\delta - \log2-(n+m)\log\alpha.
\end{equation}
Note that $e^{\Gamma_1}-1=\Lambda_1$. If $ m>100 $, then $ \frac{6}{\alpha^{2m}}< \frac{1}{2} $. Since $ |e^{\Gamma_1}-1|<1/2 $, it follows that
\begin{eqnarray}\label{Pelling12}
|\Gamma_1| < 2|e^{\Gamma_1}-1|< \dfrac{12}{\alpha^{2m}}.
\end{eqnarray}
By recalling that $ (k,n,m)=(k_i, n_i, m_i) $ for $ i=1,2 $, we get that 
\begin{eqnarray}\label{Pelling13}
\left|k_i\log\delta - \log2-(n_i+m_i)\log\alpha\right|< \dfrac{12}{\alpha^{2m_i}}
\end{eqnarray}
holds for both $ i=1, 2 $ provided $ m_3>100 $.

We apply Theorem \ref{Matveev11} on the left-hand side of \eqref{Pelling11}. First, we need to check that $ \Lambda_1 \neq 0$. Well, if it were, then 
$ \delta^{k}\alpha^{-n-m}=2$. However, this is impossible since $ \delta^{k}\alpha^{-n-m} $ is a unit while $ 2 $ is not. Thus, $ \Lambda_1 \neq 0 $, and we can apply Theorem \ref{Matveev11}. We take the data
\begin{eqnarray*}
t:=3, \quad \gamma_1:=\delta, \quad \gamma_2:=2, \quad \gamma_3:=\alpha, \quad b_1:=k, \quad b_2:=-1, \quad b_3:=-n-m.
\end{eqnarray*}
We take $ \mathbb{K}:=\mathbb{Q}(\sqrt{d}, \alpha)$ which has degree $ D\le 4 $ (it could be that $ d=5 $ in which case $ D=2 $; otherwise, $ D=4 $). Since $ \delta\ge 1+\sqrt{2}>\alpha $, the second inequality in \eqref{Pell12} tells us that $ k< n+m $, so we take $ B:=2n $. We have $ h(\gamma_1)=h(\delta)=\frac{1}{2}\log\delta $, $ h(\gamma_2)=h(2)=\log2 $ and $ h(\gamma_3)=h(\alpha)=\frac{1}{2}\log\alpha $. Thus, we can take
$ A_1:=2\log\delta $, $ A_2:=4\log2 $ and $ A_3:=2\log\alpha $. Now, Theorem \ref{Matveev11} tells us that
\begin{eqnarray*}
\log|\Lambda_1|&>&-1.4\times 30^{6}\times 3^{4.5}\times 4^{2} (1+\log 4)(1+\log (2n))(2\log\delta)(4\log 2)(2\log\alpha)\\
&>&-2.92\times 10^{13}\log\delta (1+\log(2n)).
\end{eqnarray*}
By comparing the above inequality with \eqref{Pelling11}, we get
\begin{eqnarray}
2m\log\alpha -\log 6 < 2.92\times 10^{13}\log\delta (1+\log(2n)).
\end{eqnarray}
Thus
\begin{eqnarray}\label{est}
m<6.06\times 10^{13}\log\delta (1+\log(2n)).
\end{eqnarray}
Since, $ \delta^{k}<\alpha^{n+m+6} $, we get that
\begin{eqnarray}
k\log\delta < (n+m+6)\log\alpha \le (2n+6) \log\alpha,
\end{eqnarray}
which together with the estimate \eqref{est} gives
\begin{eqnarray}
km< 5.84\times 10^{13}n(1+\log(2n)).
\end{eqnarray}
Let us record what we have proved, since this will be important later-on.

\begin{lemma}\label{DM1}
If $ x_k=L_nL_m $ and $ n\ge m $, then
$$m<6.06\times 10^{13}\log\delta (1+\log(2n)), \quad km< 5.84\times 10^{13}n(1+\log(2n)),\quad k\log\delta < 4n\log\alpha.$$
\end{lemma}

Note that we did not assume that $ m_3> 100 $ for Lemma \ref{DM1} since we have worked with the inequality \eqref{Pelling11} and not with \eqref{Pelling12}. We now again assume that $ m_3>100 $. Then the two inequalities \eqref{Pelling13} hold. We eliminate the term involving $ \log \delta $ by multiplying the inequality for $ i=1 $ with $ k_2 $ and the one for $ i=2 $ with $ k_1 $, subtract them and apply the triangle inequality as follows
\begin{eqnarray*}
&&\left|(k_2-k_1)\log 2 - (k_2(n_1+m_1)-k_1(n_2+m_2))\log\alpha\right|\\
&&= \left|k_2(k_1\log\delta -\log 2-(n_1+m_1)\log\alpha) - k_1(k_2\log\delta-\log 2 - (n_2+m_2)\log\alpha)\right|\\
&&\le k_2\left|k_1\log\delta -\log 2-(n_1+m_1)\log\alpha\right|+ k_1\left|k_2\log\delta-\log 2 - (n_2+m_2)\log\alpha\right|\\
&&\le \dfrac{12k_2}{\alpha^{2m_1}} + \dfrac{12k_1}{\alpha^{2k_2}} ~<~ \dfrac{24k_2}{\alpha^{2m_3}}.
\end{eqnarray*}
Thus,
\begin{eqnarray}\label{DM2}
|\Gamma_2|:=\left|(k_2-k_1)\log 2- (k_2(n_1+m_1)-k_1(n_2+m_2))\log\alpha\right| < \dfrac{24k_2}{\alpha^{2m_3}}.
\end{eqnarray}
We are now set to apply Theorem \ref{Matveev12} with the data
\begin{eqnarray*}
t:=2, \quad \gamma_1:=2, \quad \gamma_2:=\alpha, \quad b_1:=k_2-k_1, \quad b_2:=k_2(n_1+m_1)-k_1(n_2+m_2).
\end{eqnarray*}
The fact that $ \gamma_1=2 $ and $ \gamma_2=\alpha $ are multiplicatively independent follows because $ \alpha $ is a unit while $ 2 $ is not. We observe that $ k_2-k_1< k_2 $, whereas by the absolute value of the inequality in \eqref{DM2}, we have
\begin{eqnarray*}
\left|k_2(n_1+m_1)-k_1(n_2+m_2)\right|\le (k_2-k_1)\dfrac{\log 2}{\log\alpha}+\dfrac{24k_2}{\alpha^{2m_3}\log\alpha} <2k_2,
\end{eqnarray*}
because $ m_3> 10 $. We have that $ \mathbb{K}:=\mathbb{Q}(\alpha) $, which has $ D=2 $. So we can take
\begin{eqnarray*}
\log B_1=\max\left\{h(\gamma_1), \dfrac{|\log \gamma_1|}{2}, \dfrac{1}{2}\right\} = \log 2,
\end{eqnarray*}
and 
\begin{eqnarray*}
\log B_2=\max\left\{h(\gamma_2), \dfrac{|\log \gamma_2|}{2}, \dfrac{1}{2}\right\} = \dfrac{1}{2}.
\end{eqnarray*}
Thus,
\begin{eqnarray*}
b^{\prime}=\dfrac{|k_2-k_1|}{2\log B_2}+\dfrac{|k_2(n_1+m_1)-k_1(n_2+m_2)|}{2\log B_1} \le k_2+\dfrac{k_2}{\log 2}< 3k_2.
\end{eqnarray*}
Now Theorem \ref{Matveev12} tells us that with $$\Gamma_2= (k_2-k_1)\log 2- (k_2(n_1+m_1)-k_1(n_2+m_2))\log\alpha, $$
we have 
\begin{eqnarray*}
\log|\Gamma_2|> -24.34\times 2^{4}\left(\max\{\log(3k_2) +0.14, 10.5\}\right)^{2}\cdot (2\log 2)\cdot (1/2).
\end{eqnarray*}
Thus,
\begin{eqnarray*}
\log|\Gamma_2|>-270\left(\max\{\log(3k_2) +0.14, 10.5\}\right)^{2}.
\end{eqnarray*}
By comparing the above inequality with \eqref{DM2}, we get
\begin{eqnarray*}
2m_3\log \alpha - \log (24 k_2) <270\left(\max\{\log(3k_2) +0.14, 10.5\}\right)^{2}.
\end{eqnarray*}
If $ k_2\le 10523 $, then $ \log(3k_2) + 0.14 < 10.5 $. Thus, the last inequality above gives
\begin{eqnarray*}
2m_3\log\alpha< 270\times 10.5^{2}+\log (24\times 10523),
\end{eqnarray*}
giving $ m_3<30942 $ in this case. Otherwise, $ k_2> 10523 $, and we get
\begin{eqnarray*}
2m_3\log\alpha < 272(1+\log k_2)^{2}+\log (24k_2)< 280 (1+\log k_2)^{2},
\end{eqnarray*}
which gives
\begin{eqnarray*}
m_3< 160 (1+\log k_2)^{2}.
\end{eqnarray*}
We record what we have proved
\begin{lemma}\label{DM22}
If $ m_3>100 $, then either 
\begin{itemize}
\item[(i)] $ k_2\le 10523 $ and $ m_3< 30942 $ or 
\item[(ii)] $ k_2> 10523 $, in which case $ m_3< 160(1+\log k_2)^{2} $.
\end{itemize}
\end{lemma}

Now suppose that some $ m $ is fixed in \eqref{Problem}, or at least we have some good upper bounds on it. We rewrite \eqref{Problem} using \eqref{Binet} and \eqref{Pellxterm} as
\begin{eqnarray*}
\dfrac{1}{2}(\delta^{k}+\eta^{k})= L_m(\alpha^{n}+\beta^{n}),
\end{eqnarray*}
so
\begin{eqnarray*}
\delta^{k}\left(2L_m\right)^{-1}\alpha^{-n}-1 =-\dfrac{1}{2L_m}\eta^{k}\alpha^{-n}+(\beta\alpha^{-1})^{n}.
\end{eqnarray*}
Since $ m\ge 1 $, $ \beta =-\alpha^{-1} $, $ \eta = \varepsilon \delta^{-1} $ and $ \delta^{k}> \alpha^{n+m-1} $, we get
\begin{eqnarray*}
\left|\delta^{k}\left(2L_m\right)^{-1}\alpha^{-n}-1\right|&\le& \dfrac{1}{2L_m\delta^{k}\alpha^{n}}+\dfrac{1}{\alpha^{2n}}\quad \le\quad\dfrac{\alpha}{\alpha^{2(n+m)}}+\dfrac{1}{\alpha^{2n}}\\
&\le& \dfrac{\alpha+1}{\alpha^{2n}} \quad < \quad \dfrac{6}{\alpha^{2n}},
\end{eqnarray*}
where we have used the fact that $ n\ge m\ge 0 $ and $ \alpha+1<6 $. Hence,
\begin{eqnarray}\label{DM3}
|\Lambda_3|:= \left|\delta^{k}\left(2L_m\right)^{-1}\alpha^{-n}-1\right|< \dfrac{6}{\alpha^{2n}}.
\end{eqnarray}
We assume that $ n_3>100 $. In particular, $ \frac{6}{\alpha^{2n}}< \frac{1}{2} $ for $ n\in\{n_1, n_2\} $, so we get by the previous argument that
\begin{eqnarray}\label{DM4}
|\Gamma_3|:=\left|k\log\delta -\log(2L_m)-n\log\alpha\right| < \dfrac{12}{\alpha^{2n}}.
\end{eqnarray}
We are now set to apply Theorem \ref{Matveev11} on the left-hand side of \eqref{DM3} with the data
\begin{eqnarray*}
t:=3, \quad \gamma_1:=\delta, \quad \gamma_2:=2L_m, \quad \gamma_3:=\alpha, \quad b_1:=k, \quad b_2:= -1, \quad b_3: = -n.
\end{eqnarray*} 
First, we need to check that $ \Lambda_3: = \delta^{k}(2L_m)^{-1}\alpha^{-n}-1 \neq 0 $. If not, then $ \delta^{k}=2L_m\alpha^{m} $. The left-hand side belongs to the field $ \mathbb{Q}(\sqrt{d}) $ but not rational while the right-hand side belongs to the field $ \mathbb{Q}(\sqrt{5}) $. This is not possible unless $ d=5 $. In this last case, $ \delta $ is a unit in $ \mathbb{Q}(\sqrt{5}) $ while $ 2L_m$ is not a unit in $ \mathbb{Q}(\sqrt{5}) $ since the norm of this first element is $ 4L_m^2\neq \pm 1 $. So, $ \Lambda_3\neq 0 $. Thus, we can apply Theorem \ref{Matveev11}. We have the field $ \mathbb{K}:=\mathbb{Q}(\sqrt{d}, \sqrt{5})$ which has degree $ D\le 4 $. We also have 
\begin{eqnarray*}
h(\gamma_2)&=&h(2L_m) = h(2)+h(L_m) \\ &\le& \log 2 + (m+1)\log\alpha < 2+m\log\alpha\\
&\le&2.92\times 10^{13}\log\delta(1+\log (2n)) \quad \text{by \eqref{est}}.
\end{eqnarray*}
So, we take
\begin{eqnarray*}
h(\gamma_1) = \dfrac{1}{2}\log\delta, \quad h(\gamma_2)  = 2.92\times 10^{13}\log\delta(1+\log (2n)) \quad \text{and} \quad h(\gamma_3) = \dfrac{1}{2}\log\alpha.
\end{eqnarray*}
Then,
\begin{eqnarray*}
A_1: = 2\log\delta, \quad A_2:=1.18\times 10^{14}\log\delta(1+\log (2n)) \quad \text{and} \quad A_3:=2\log \alpha.
\end{eqnarray*}
Then, by Theorem \ref{Matveev11} we get
\begin{eqnarray*}
\log|\Lambda_3|&>& -1.4\times 30^{6}\times 3^{4.5}\times 4^{2}(1+\log 4)(1+\log n)(2\log\delta)\\&& \times(1.18\times 10^{14}\log\delta(1+\log (2n)))(2\log\alpha)\\
&>&-8.6\times 10^{26}(1+\log (2n))^2(\log\delta)^2\log\alpha.
\end{eqnarray*}
Comparing the above inequality with \eqref{DM3}, we get
\begin{eqnarray*}
2n\log \alpha -\log 6 < 8.6\times 10^{26}(1+\log (2n))^2(\log\delta)^2\log\alpha,
\end{eqnarray*}
which implies that
\begin{eqnarray}
n< 4.3\times 10^{26}(1+\log (2n))^2(\log\delta)^2.
\end{eqnarray}
We record what we have proved.
\begin{lemma}\label{DM11}
If $ x_k=L_nL_m $ with $ n\ge m\ge 1 $, then we have
\begin{eqnarray*}
n< 4.3\times 10^{26}(1+\log (2n))^2(\log\delta)^2.
\end{eqnarray*}
\end{lemma}
Note that we did not use the assumption that $ m_3>100 $ of that $ n_3>100 $ for Lemma \ref{DM11} since we worked with the inequality \eqref{DM3} not with the inequality \eqref{DM4}. We now assume that $ n_3>100 $ and in particular \eqref{DM4} holds for $ (k, n, m)= (k_i, n_i, m_i) $ for both $ i=1,2 $. By the previous procedure, we also eliminate the term involving $ \log\delta $ as follows
\begin{eqnarray}\label{lll1}
&\left| k_2\log (2L_{m_1})-k_1\log (2L_{m_2}) - (k_2n_1-k_1n_2)\log\alpha\right| < \dfrac{12k_2}{\alpha^{2n_1}} + \dfrac{12 k_1}{\alpha^{2n_2}}
<\dfrac{24k_2}{\alpha^{2n_3}}.
\end{eqnarray}
We assume that $ \alpha^{2n_3}>48k_2 $. If we put
\begin{eqnarray*}
\Gamma_4: =  k_2\log (2L_{m_1})-k_1\log (2L_{m_2}) - (k_2n_1-k_1n_2)\log\alpha,
\end{eqnarray*}
we have that $ |\Gamma_4|<1/2 $. We then get that
\begin{eqnarray}\label{DMA2}
|\Lambda_4|:=|e^{\Gamma_4}-1|< 2|\Gamma_4|< \dfrac{48k_2}{\alpha^{2n_3}}.
\end{eqnarray}
We apply Theorem \ref{Matveev11} to
\begin{eqnarray*}
\Lambda_4:=(2L_{m_1})^{k_2}(2L_{m_2})^{-k_1}\alpha^{-(k_2n_1-k_1n_2)}-1.
\end{eqnarray*}
First, we need to check that $ \Lambda_4\neq 0 $. Well, if it were, then it would follow that 
\begin{eqnarray}\label{musua}
\dfrac{L_{m_1}^{k_2}}{L_{m_2}^{k_1}}=2^{k_1-k_2}\alpha^{k_2n_1-k_1n_2}.
\end{eqnarray}
We consider the following Lemma.
\begin{lemma}\label{musu1}
The equation \eqref{musua} has only many small positive integer solutions $(k_i, n_i, m_i)$ for $ i=\{1,2\} $ with  $ k_1<k_2 $ and $  m_1 \le m_2 \le 6 $. Futhermore, none of these solutions lead to a valid solution to the original Diophantine equation \eqref{Problem}.
\end{lemma}
\begin{proof}
$  $We suppose that \eqref{musua} holds and assume that $ \gcd (k_1, k_2)=1 $. Since $ \alpha^{k_2n_1-k_1n_2} \in\mathbb{Q} $, it follows  $ k_2n_1=k_1n_2 $. Thus, if one of the $ n_1, ~ n_2 $ is zero, so is the other. Since $ n_i\ge m_i $ for $ i\in\{1,2\} $, it follows that $ n_1=n_2=0 $, $ m_1=m_2=0 $, so $ x_{k_1}=x_{k_2} $, therefore $ k_1=k_2 $ a contradiction. Thus, $ n_1 $ and $ n_2 $ are both positive integers. Next $ {L_{m_1}^{k_2}}/{L_{m_2}^{k_1}}=2^{k_1-k_2}<1 $. Thus, $ {L_{m_1}^{k_2}}<{L_{m_2}^{k_1}}< {L_{m_2}^{k_2}}$, so $ L_{m_1}< L_{m_2} $. This implies that either $ (m_1, m_2)=(1,0) $ or $ m_1<m_2 $. The case $ (m_1, m_2)=(1,0) $ gives 
$ 1/2^{k_1}=2^{k_1-k_2} $. Thus, $ k_2=2k_1 $ and since
 $ \gcd(k_1, k_2)=1 $, we get 
 $ k_1=1, ~ k_2 =2 $, so  $ n_2=2{n_1} $. But then 
 $ x_2=x_{k_2}=L_{n_2}L_{m_2} = L_{2n_1}L_0=2L_{2{n_1}}$ is even, a contradiction since $ x_2=2x_1 \pm 1 $ (by Example \ref{exa:2} (i)) is odd. Thus, $ m_1< m_2 $. If $ m_2>6 $, the Carmichael Primitive Divisor Theorem for Lucas numbers shows that $ L_{m_2} $ is divisible by a prime $ p>7 $ which does not divide $ L_{m_1} $. This is impossible since it contradicts the assumption that \eqref{musua} holds. Thus, $ m_2 \le 6 $. Further since $ {L_{m_1}^{k_2}}/{L_{m_2}^{k_1}}=1/2^{k_2-k_1} $ it follows that $ {L_{m_1}^{k_1}}\mid {L_{m_1}^{k_2}}\mid {L_{m_2}^{k_1}}$, so $ L_{m_1}\mid L_{m_2} $. So, there are three cases that we analyse:
 
 \medskip
 
 \textbf{Case 1}. $ m_1=0 $, $ m_2\in\{3, 6\} $. If $ (m_1, m_2)=(0,3) $, then $ 2^{k_2}/4^{k_1}= 1/2^{2k_1-k_2}= 1/2^{k_2-k_1} $. This gives $ 2k_2=3k_1 $ and since $ k_1 $ and $ k_2 $ are coprime, it follows that $ k_1=2 $ and $ k_2=3 $. Then 
 $ x_2=x_{k_1}=L_{n_1}L_{m_1} = L_{n_1}L_0=2L_{{n_1}}$ is even, a contradiction since $ x_2=2x_1 \pm 1 $ is odd. If $ (m_1, m_2)=(0,6) $, then $ 2^{k_2}/18^{k_1} =1/2^{k_2-k_1}$, which is impossible  since by looking at the exponent of $ 3 $ we would get $ k_1=0 $, a contradiction.
 
 \medskip
 
 \textbf{Case 2}. $ m_1=2 $ and $ L_{m_2} $ is a power of $ 2 $. The case $ m_2=0 $ has been treated so the only other case left is $ m_2=3 $. In this case, $ 1/4^{k_1}=1/2^{k_2-k_1} $, giving $ k_2=3k_1 $. Thus, since $ \gcd(k_1, k_2)=1 $, then $ k_1=1 $ and $ k_2=3 $. Since $ k_2n_1=k_1n_2 $, we get $ n_2=3n_1 $. Thus, $ x_1=L_{n_1}L_{1} = L_{n_1}$ and $ x_3=L_{3n_1}L_3=4L_{3n_1} $. Now $ x_3=x_1(4x_1^2\pm 3) $ (by Example \ref{exa:2} (ii)) and the second factor is odd, so the power of $ 2 $ dividing $ 4L_{3n_1} $ divides $ x_1=L_{n_1} $. But $ 4L_{3n_1} $ is a multiple of $ 8 $ since $ L_{3n_1} $ is even. Thus, $ 8\mid L_{n_1} $, which is false.
 
 \medskip
 
 \textbf{Case 3}. $ m_1=2 $ and $ m_2=6 $. We get $ 3^{k_2}/(2.3^2)^{k_1}=1/2^{k_2-k_1} $. Looking at the exponent of $ 3 $, we get $ k_2=2k_1 $ and loking at the exponent of $ 2 $ we also get $ k_2=2k_1 $, so $ k_1=1 $ and $ k_2=2 $. Also, $ n_2=2n_1 $. Thus, $ x_1=L_{n_1}L_{m_1}=3L_{n_1} $ and $ x_2=L_{n_2}L_{m_2}=18L_{2n_1} $ is even, a contradiction with the fact that $ x_2=2x_1^2\pm 1 $ is odd.
\end{proof}
 So, by Lemma \ref{musu1} we have $ \Lambda_4 \neq 0$. Thus, we can now apply Theorem \ref{Matveev11} with the data
\begin{eqnarray*}
&t:=3, \quad \gamma_1:=2L_{m_1}, \quad \gamma_2:= 2L_{m_2}, \quad \gamma_3:=\alpha, \quad b_1=k_2, \\ & b_2:=-k_1, \quad b_3:=-(k_2n_1-k_1n_2).
\end{eqnarray*}
We have $ \mathbb{K}:=\mathbb{Q}(\sqrt{5}) $ which has degree $ D:=2 $. Also, using \eqref{DMA1}, we can take $ B:=4n_4^{2} $. We can also take $ A_1:=2(2+m_1\log \alpha) \le 4m_1\log\alpha $, $ A_2:=2(2+m_2\log \alpha) \le 4m_2\log\alpha $ and $ A_3:=\log\alpha $. Theorem \ref{Matveev11} gives that
\begin{eqnarray*}
\log |\Lambda_4|&>&-1.4 \times 30^{6}\times 3^{4.5}\times 2^{2}(1+\log 2)(1+\log (4n_4^2))(4m_1\log \alpha)(4m_2\log \alpha)\log\alpha,\\
&>&-3.44\times 10^{12}m_1m_2(1+\log (2n_4)).
\end{eqnarray*}
By comparing this with the inequality \eqref{DMA2}, we get
\begin{eqnarray*}
2n_3\log\alpha - \log (48k_2) < 3.44\times 10^{12}m_1m_2(1+\log (2n_4)).
\end{eqnarray*}
Since $ k_2<4n_4 $ and $ n_4>10 $, we get that $ \log (48 k_2) <2(1+\log (2n_4)) $. Thus,
\begin{eqnarray}\label{DMA3}
n_3< 3.58\times 10^{12}m_1m_2(1+\log (2n_4)).
\end{eqnarray}
All this was done under the assumption that $ \alpha^{2n_3}>48k_2 $. But if that inequality fails, then
\begin{eqnarray*}
n_3<c_1 \log(48k_2) < 12(1+\log (2n_4)),
\end{eqnarray*}
which is much better than \eqref{DMA3}. Thus, \eqref{DMA3} holds in all cases. Next, we record what we have proved.
\begin{lemma}\label{DM33}
Assuming that $ n_3>100 $, then we have
\begin{eqnarray*}
n_3<3.58\times 10^{12}m_1m_2(1+\log (2n_4)).
\end{eqnarray*}
\end{lemma}
We now start finding effective bounds for our variables.

\textbf{ Case 1.} $ m_4\le 100 $.

Then $ m_1<100 $ and $ m_2<100 $. By Lemma \ref{DM33}, we get that
\begin{eqnarray*}
n_3< 3.58\times 10^{16}(1+\log (2n_4)).
\end{eqnarray*}
By Lemma \ref{DM1}, we get
\begin{eqnarray*}
\log\delta < 4n_3\log\alpha < 6.89\times 10^{16}(1+\log (2n_4)).
\end{eqnarray*}
By the inequality \eqref{Pell12}, we have that
\begin{eqnarray*}
n_4&\le& n_4+m_4-1\\
&< & k_2c_1\log\delta\\
&<& 1.72\times 10^{27}c_1(1+\log(2n_4))^2(\log\delta)^3 \quad \text{(by \eqref{DMA1} and Lemma \ref{DM11})}\\ 
&<&\dfrac{1}{\log\alpha}(1.72\times 10^{27}(1+\log(2n_4))^2)(6.89\times 10^{16}(1+\log (2n_4)))^3\\
&<& 1.17\times 10^{78}\log (1+\log(2n_4))^5.
\end{eqnarray*}
With the help of \textit{Mathematica}, we get that $ n_4< 4.6\times 10^{89} $. Thus, using \eqref{DMA1}, we get
\begin{eqnarray*}
\max\{k_2, n_4\} < 4.6\times 10^{89}.
\end{eqnarray*}
We record what we have proved.
\begin{lemma}\label{DM44}
If $ m_4:=\max\{m_1, m_2\} \le 100 $, then
\begin{eqnarray*}
\max\{k_2, n_4\} < 4.6\times 10^{89}.
\end{eqnarray*}
\end{lemma}

\textbf{ Case 2.} $ m_4> 100 $.

Note that either $ m_3\le 100 $ or $ m_3>100 $
case in which by Lemma \ref{DM22} and the inequality \eqref{DMA1}, we have $ m_3\le 160 (1+\log(4n_4))^{2}  $  provided that $ m_4> 10000 $, which we now assume.

We let $ i\in\{1,2\} $ be such that $ m_i=m_3 $ and $ j $ be such that $ \{i,j\}=\{1,2\} $. We assume that $ n_3> 100 $. We work with \eqref{DM4} for $ i $ and \eqref{Pelling13} for $ j $ and noting the conditions $ n_i>100 $ and $ m_j=m_4>100 $ are fullfilled. That is,
\begin{eqnarray*}
\left|k_i\log\delta +\log(2L_{m_i})-n_i\log\alpha\right| &<& \dfrac{12}{\alpha^{2n_i}},\\
\left|k_j\log\delta - \log2-(n_j+m_j)\log\alpha\right|&<& \dfrac{12}{\alpha^{2m_j}}.
\end{eqnarray*}
By a similar procedure as before, we eliminate the term involving $ \log\delta $.  We multiply the first inequality by $ k_j $, the second inequality by $ k_i $, subtract the resulting inequalities and apply the triangle inequalty to get
\begin{eqnarray}\label{MAA}
\left|k_j\log(2L_{m_i})-k_i\log 2-(k_jn_i-k_i(n_j+m_j))\log\alpha\right|&<& \dfrac{12k_j}{\alpha^{2m_i}}+\dfrac{12k_i}{\alpha^{2l_j}}\nonumber\\
&<& \dfrac{24k_2}{\alpha^{2\min\{n_i, m_j\}}}.
\end{eqnarray}
Assume that $ \alpha^{2\min\{n_i, m_j\}} > 48k_2 $. We put
\begin{eqnarray*}
\Gamma_5:= k_j\log(2L_{m_i})-k_i\log 2-(k_jn_i-k_i(n_j+m_j))\log\alpha.
\end{eqnarray*}
We can write $ \Lambda_5:=(2L_{m_i})^{k_j} 2^{-k_i}\alpha^{(k_jn_i-k_i(n_j+m_j))}-1$. Under the above assumption and using \eqref{MAA}, we get that
\begin{eqnarray}\label{MAA1}
|\Lambda_5|= |e^{\Gamma_5}-1| < 2|\Gamma_5| < \dfrac{48k_2}{\alpha^{2\min\{n_i, m_j\}}}.
\end{eqnarray}

We are now set to apply Theorem \ref{Matveev11} on $ \Lambda_5 $. First, we need to check that $ \Lambda_5\neq 0 $. Well, if it were, then we would get that
\begin{eqnarray}\label{MA}
L_{m_i}^{k_j}=2^{k_i-k_j}\alpha^{(k_jn_i-k_i(n_j+m_j))}.
\end{eqnarray}
We consider the following lemma.
\begin{lemma}\label{musu2}
The equation \eqref{MA} has only many small positive integer solutions \\$(k_i, k_j, n_i, n_j, m_i, m_j)$ for $ i,j=\{1,2\} $ with  $ k_1<k_2 $ and $  m_1 \le m_2 \le 6 $. Futhermore, none of these solutions lead to a valid solution to the original Diophantine equation \eqref{Problem}.
\begin{proof}
Suppose that \eqref{MA} holds and assume that $ \gcd(k_1, k_2)=1 $. Since $ \alpha^{(k_jn_i-k_i(n_j+m_j))} \in \mathbb{Q} $, then $ k_jn_i=k_i(n_j+m_j) $. Next $ L_{m_i}^{k_j}=2^{k_i-k_j} $. Thus, $ k_i \ge k_j $, so $ i=2, ~j=1 $, $ k_2>k_1 $ and $ m_2 \neq 1 $. Since $ L_{m_2}>1 $ is a power of $ 2 $, it follows that $ m_2 \in \{0,3\} $. Suppose $ m_2=0 $. Then $ L_{m_2}^{k_1}=2^{k_1}=2^{k_2-k_1} $, so $ k_2=2k_1 $. Hence, $ k_1=1 $ and $ k_2=2 $. Further, $ n_2=2(n_1+m_1) $. Thus, $ x_2=x_{k_2}=L_{n_2}L_{m_2}=2L_{2(n_1+m_1)} $ is even, which false because $ x_2=2x_1^2\pm 1 $ is odd. Suppose next that $ m_2=3 $. Then $ 4^{k_1}= 2^{k_2-k_1} $. Thus, $ k_2=3k_1 $, so $ k_1=1  $ and $ k_2 =3 $. Next, $ n_2 = 3 (n_1+m_1) $. Hence,  $ x_1=x_{k_1}=L_{n_1}L_{m_1}$ and   $ x_3=x_{k_2}=L_{n_2}L_{m_2}=4L_{3(n_1+m_1)} $. By the previous argument in the proof of Lemma \ref{musu1}, $ 8 $ divides $ x_3=x_1(4x_1^2\pm 1) $, so $ 8\mid x_1 $. Since $ x_1 = L_{n_1}L_{m_1} $ and $ 8 \nmid L_n$ for any $ n $, it follows that $ L_{n_1} $ and $ L_{m_1} $ are both even. Thus, $ 3\mid n_1 $, $ 3\mid m_1 $. Further, one of $ L_{n_1} $, $ L_{m_1} $ is a multiple of $ 4 $, so one of $ n_1, ~m_1 $ is odd. Suppose both are odd. Then $ 4 \mid L_{n_1}$, $ 4 \mid L_{m_1}$ so $ 16\mid x_1\mid x_3\mid 4L_{3(n_1+m_1)} $. This implies that $ 4\mid L_{3(n_1+m_1)} $, which is false because $ 3(n_1+m_1) $ is an even multiple of $ 3 $, and $ 2 \Vert L_{6m} $ for any $ m $. Suppose now that one of $ n_1, ~m_1 $ is an even multiple of $ 3 $, and the other is odd. Then $ \text{ord}_{2}(x_1)=3 $, where $ \text{ord}_{2}(x) $ is the exponent at which $ 2 $ appears in the factorization of $ x $. Hence,
$$
3=\text{ord}_2(x_3)=\text{ord}_2(4L_{3(n_1+m_1)})=2+\text{ord}_2(L_{3(n_1+m_1)}),
$$
giving $ \text{ord}_2(L_{3(n_1+m_1)})=1 $, which is again false since $ 3(n_1+m_1) $ is an odd multiple $ 3 $, so a number of the form $ 3+6m $, and for such numbers we have $ 4 \Vert L_{3+6m} $. Hence, in all instances we have gotten a contradiction.
\end{proof}
\end{lemma}
Thus, by Lemma \ref{musu2} we have that1 $ \Lambda_5\neq 0 $. So, we can apply Theorem \ref{Matveev11} with the data
\begin{eqnarray*}
&t:=3, \quad \gamma_1:=2L_{m_i}, \quad \gamma_2:= 2\quad \gamma_3:=\alpha\quad 
b_1:=k_j, \\ &b_2:= -k_i, \quad b_3:= -(k_jn_i-k_i(n_j+m_j)).
\end{eqnarray*}
From the previous calculations, we know that $ \mathbb{K}:=\mathbb{Q}(\sqrt{2} )$ which has degree $ D=2 $ and $ A_1:= 4m_i\log \alpha $, $ A_2:=2\log 2 $ and $ A_3:=\log\alpha $.
We also take $ B:=4n_4^2 $. By Theorem \ref{Matveev11}, we get that
\begin{eqnarray*}
\log |\Lambda_5|&>&-1.4 \times 30^{6}\times 3^{4.5}\times 2^{2}(1+\log 2)(1+\log (4n_4^2))(4m_i\log \alpha)(2\log 2)\log\alpha,\\
&>&-5.18\times 10^{12}m_i(1+\log (2n_4)).
\end{eqnarray*}
Comparing the above inequality with \eqref{MAA1}, we get
\begin{eqnarray*}
2\min\{n_i, m_j\}\log\alpha -\log(48k_2) < 5.12\times 10^{12}m_i(1+\log (2n_4)).
\end{eqnarray*}
Since $ m_4>100 $, we get using \eqref{DMA1} ( $ k_2<4n_4 $) that,
\begin{eqnarray*}
\min\{n_i, n_j\}<5.38\times 10^{12}(160(1+\log(4n_4))^{2})(1+\log (2n_4))+ \dfrac{c_1}{2}\log (192n_4),
\end{eqnarray*}
which implies that
\begin{eqnarray}\label{DMAA1}
\min\{n_i, m_j\}<1.72\times 10^{15}(1+\log (2n_4))^3.
\end{eqnarray}
All this was under the assumptions that $ n_4>10000 $, and that $ \alpha^{2\min\{n_i, m_j\}} > 48k_2 $. But, still under the condition that $ n_4>10000 $, if  $ \alpha^{2\min\{n_i, m_j\}} < 48k_2 $, then we get an inequality for $ \min\{n_i, n_j\} $ which is even much better than \eqref{DMAA1}. So, \eqref{DMAA1} holds provided that $ n_4>10000 $. Suppose say that $ \min\{n_i, m_j\}= m_j $. Then we get that
\begin{eqnarray*}
m_3<160(1+\log(4n_4))^{2}, \quad m_4< 1.72\times 10^{15}(1+\log (2n_4))^3.
\end{eqnarray*}
By Lemma \ref{DM33}, since $ m_3>100 $, we get
\begin{eqnarray*}
n_3&<&(3.58\times 10^{12})(160(1+\log(4n_4))^{2})(1+\log (2n_4))\\
&&\times 1.72\times 10^{15}(1+\log (2n_4))^3\\
&<&1.98\times 10^{30}(1+\log (2n_4))^{6}.
\end{eqnarray*}
Together with Lemma \ref{DM1}, we get
\begin{eqnarray*}
\log\delta < 3.80\times 10^{30}(1+\log (2n_4))^{6},
\end{eqnarray*}
which together with Lemma \ref{DM11} gives
\begin{eqnarray*}
n_4&<& 4.30\times 10^{26}(1+\log (2n_4))^{2}(3.80\times 10^{30}(1+\log (2n_4))^{6})^2,
\end{eqnarray*}
which implies that
\begin{eqnarray}
n_4< 6.21\times 10^{87}(1+\log (2n_4))^{14}.
\end{eqnarray}
With the help of \textit{Mathematica} we get that $ n_4< 1.30\times 10^{122} $. This was proved under the assumption that $ n_4>10000 $, but the situation $ n_4\le 10000 $ already provides a better bound than  $ n_4< 1.30\times 10^{122} $. Hence,
\begin{eqnarray}
\max\{k_2, n_1, n_2\}< 1.30\times 10^{122}.
\end{eqnarray}

This was when $ m_j=\min\{n_i, m_j\} $. Now we assume that $ n_i=\min\{n_i, m_j\} $. Then we get
\begin{eqnarray*}
n_i< 1.72\times 10^{15}(1+\log (2n_4))^3.
\end{eqnarray*}
By Lemma \ref{DM1}, we get that
\begin{eqnarray*}
\log\delta < 3.31\times 10^{15}(1+\log (2n_4))^3.
\end{eqnarray*}
Now by Lemma \ref{DM11} together with Lemma \ref{DM1} to bound $ l_4 $ give
\begin{eqnarray*}
n_4&<& 4.30\times 10^{26}(1+\log(2n_4)))^{2}(3.31\times 10^{15}(1+\log (2n_4))^3)^2\\
&<&4.72\times 10^{57}(1+\log (2n_4))^{10}.
\end{eqnarray*}
This gives, $ n_4< 2.44\times 10^{80} $ which is a better bound than $ 1.30\times 10^{122} $. We record what we have proved.
\begin{lemma}\label{DM55}
If $ m_4: = \max\{m_1, m_2\} > 100 $ and $ n_3:=\min\{n_1, n_2\} > 100 $, then
\begin{eqnarray*}
\max\{k_2, n_1, n_2\} < 1.30\times 10^{122}.
\end{eqnarray*}
\end{lemma}

It now remains the case when $ m_4> 100 $ and $ n_3\le 100 $. But then, by Lemma \ref{DM1}, we get $ \log\delta <192 $ and now Lemma \ref{DM1} together with Lemma \ref{DM11} give
\begin{eqnarray*}
n_4<1.56\times 10^{31}(1+\log (2n_4))^2,
\end{eqnarray*}
which implies that $ n_4< 10^{36} $ and further $ \max\{k_1, n_1, n_2 \} < 10^{40} $.
We record what we have proved.
\begin{lemma}\label{DM66}
If $ m_4 > 100 $ and $ n_3 \le 100 $, then
\begin{eqnarray*}
\max\{k_1, n_1, n_2\} < 10^{40}.
\end{eqnarray*}
\end{lemma}

\section{The final computations}
\subsection{The first reduction}\label{RED}
In this subsection we reduce the bounds for $ k_1, ~ m_1, ~ n_1 $ and $k_2, m_2, ~ n_2 $ to cases that can be computationally treated. For this we return to the inequalities for $ \Gamma_2, ~\Gamma_4 $ and $ \Gamma_5 $.

We return to \eqref{DM2} and we set $ s:= k_2-k_1 $ and $ r:=k_2(n_1+m_1)-k_1(n_2+m_2) $ and divide both sides by $ s\log\alpha $ to get
\begin{eqnarray}\label{DMAA12}
\left|\dfrac{\log 2}{\log\alpha} - \dfrac{r}{s}\right|< \dfrac{24k_2}{\alpha^{2m_3}s\log\alpha}.
\end{eqnarray}
We assume that $ l_3 $ is so large that the right-hand side of the inequality in  \eqref{DMAA12} is smaller than $ 1/(2s^2) $. This certainly holds if
\begin{eqnarray}\label{DNA1}
\alpha^{2m_3}>48k_2^2/\log\alpha.
\end{eqnarray}
Since $ k_2< 1.3\times 10^{122} $, it follows that the last inequality \eqref{DNA1} holds provided that $ m_3\ge 589 $, which we now assume. In this case $ r/s $ is a convergent of the continued fraction of $ \tau: = \frac{ \log 2 }{\log \alpha} $ and $ s< 1.30\times 10^{122} $. We are now set to apply Lemma \ref{Legendre}.

We write $ \tau: =[a_0; a_1, a_2, a_3, \ldots] = [1, 2, 3, 1, 2, 3, 2, 4, 2, 1, 2, 11, 2, 1, 11, 1, 1, 134, 2, 2, \ldots] $ for the continued fraction of $ \tau $ and $ p_k/q_k $ for the $ k- $th convergent. We get that $ r/s=p_j/q_j $ for some $ j\le 237  $. Furthermore, putting $ a(M):=\max\{a_j: j=0,1, \ldots, 237\} $, we get $ a(M):=880 $. By Lemma \ref{Legendre}, we get
\begin{eqnarray*}
\dfrac{1}{882s^2}=\dfrac{1}{(a(M)+2)s^{2}}\le \left|\tau-\dfrac{r}{s}\right|<  \dfrac{24k_2}{\alpha^{2m_3}s\log\alpha},
\end{eqnarray*}
giving
\begin{eqnarray*}
\alpha^{2m_3}< \dfrac{882\times 24 k_2^2}{\log\alpha} < \dfrac{882\times 24 \times (1.30\times 10^{122})^2}{\log\alpha},
\end{eqnarray*}
leading to $ m_3\le 1190 $. We record what we have just proved.
\begin{lemma}\label{DM77}
We have $ m_3:=\min\{m_1, m_2\}\le 1190 $.
\end{lemma}

If $ m_1=m_3 $, then we have $ i=1 $ and $ j=2 $, otherwise $ m_2=m_3 $ implying that we have $ i=2 $ and $ j=1 $. In both cases, the next step is the application of Lemma \ref{LLL}  (LLL algorithm) for \eqref{MAA}, where $ n_i < 1.30\times 10^{112} $ and $ |k_jn_i-k_i(n_j +m_j)|<10^{116} $. For each $ m_j\in [1, 1190] $ and 
\begin{eqnarray}\label{DAM1}
\Gamma_5:=k_j\log(2L_{m_i})-k_i\log 2-(k_jn_i-k_i(n_j+m_j))\log\alpha,
\end{eqnarray}
we apply the LLL-algorithm on $ \Gamma_3 $ with the data
\begin{eqnarray*}
t:=3, \quad \tau_1:=\log(2L_{m_i}), \quad \tau_2:= \log 2,\quad  \tau_3:= \log\alpha\\
x_1:= k_j, \quad x_2:= -k_i, \quad x_3:= k_jn_i-k_i(n_j+m_j).
\end{eqnarray*}
Further, we set $ X:=10^{116} $ as an upper bound to $ |x_i| $ for $ i=1,2,3 $, and $ C:=(5X)^5 $. A computer search in \textit{Mathematica} allows us to conclude, together with the inequality \eqref{MAA}, that
\begin{eqnarray}
2\times 10^{-480}< \min_{1\le \min\{n_i,m_j\}\le 1190}|\Gamma_5|< \dfrac{24k_2}{\alpha^{2\min\{n_i, m_j\}}}.
\end{eqnarray}
Thus, $ \min\{n_i, m_j\}\le 1419 $. We assume first that $ i=1 , ~j=2$. Thus, $ n_1 \le 1419 $ or $ m_j=\min\{n_i, m_j\} \le 1419 $.

Next, we suppose that $ m_j=\min\{n_i, m_j\} \le 1419 $. Since $ m_1:=m_3\le 1190 $, we have
\begin{eqnarray*}
m_3:=\min\{m_1, m_2\} \le 1190  \quad \text{and} \quad m_4: = \max\{m_1, m_2\}\le 1419.
\end{eqnarray*}
Now, returning to the inequality \eqref{lll1} which involves
\begin{eqnarray}\label{DAM2}
\Gamma_4: =  k_2\log (2L_{m_1})-k_1\log (2L_{m_2}) - (k_2n_1-k_1n_2)\log\alpha, 
\end{eqnarray}
we use again the LLL algorithm to estimate the lower bound for $ |\Gamma_4| $ and thus, find a bound for $ n_1 $ that is better than the one given in Lemma \ref{DM55}. We distinguish the cases $ m_3< m_4 $ and $ m_3=m_4 $.
\subsubsection{The case $ m_3< m_4 $}
We take $ m_1:=m_3\in[1, 1190] $ and $ m_2:=m_4\in[m_3+1, 1419] $ and apply Lemma \ref{LLL} with the data:
\begin{eqnarray*}
&t:=3, \quad \tau_1:=2L_{m_1}, \quad \tau_2:=2L_{m_2}, \quad \tau_3:=\log\alpha,\\  &x_1:=k_2, \quad x_2:=-k_1, \quad x_2:=k_1n_2-k_2n_1.
\end{eqnarray*}
We also put $ X:=10^{116} $ and $ C:=(20X)^{9} $. After a computer search in \textit{Mathematica} together with the inequality \eqref{lll1}, we can confirm that
\begin{eqnarray*}
2\times 10^{-1120}\le \min_{\substack{1\le m_3\le 1190\\m_3+1\le m_4\le 1419 }} |\Gamma_4|< 24k_2\alpha^{-2n_3}.
\end{eqnarray*}
This leads to the inequality 
\begin{eqnarray*}
\alpha^{2n_3}< 12\times 10^{1120}k_2.
\end{eqnarray*}
Sustituting for the bound $ k_2 $ given in Lemma \ref{DM55}, we get that $ n_1:=n_3 \le 2950 $.
\subsubsection{The case $ m_3=m_4 $}. In this case $ m_1=m_2\le 1419 $ and we have
\begin{eqnarray*}
\Gamma_4: =  (k_2-k_1)\log (2L_{m_1}) - (k_2n_1-k_1n_2)\log\alpha \neq 0.
\end{eqnarray*}
This is similar to the case we have handled in the previous steps and yields the bound on $ n_1 $ which is less than $ 2950 $. So in both cases we have $ n_1\le 2950 $. From the fact that
\begin{eqnarray*}
\log\delta \le k_1\log\delta \le 4n_1\log\alpha < 5678,
\end{eqnarray*}
and by considering the inequality given in Lemma \ref{DM11}, we conclude that $$ n_2 < 1.4\times 10^{34}(1+\log(2n_2))^2, $$ which with the help of \textit{Mathematica} yields $ n_2< 1.12\times 10^{38} $. We summarise the first cycle of our reductions.
\begin{eqnarray}\label{pil}
\max\{k_1, m_1\}\le n_1< 2950 \quad \text{and} \quad \max\{k_2, m_2\}\le n_2< 1.12\times 10^{38}.
\end{eqnarray}
From \eqref{pil}, we note that the upper bound on $ n_2 $ represents a very good reduction of the bound given in Lemma \ref{DM55}. Hence, we expect that if we restart our reduction cycle with the new bound on $ n_2 $, then we get better bounds on $ n_1 $ and $ n_2 $. Thus, we return to the inequality \eqref{DMAA12} and take $ M:=1.12\times 10^{38} $. A computer seach in \textit{Mathematica} reveals that
\begin{eqnarray*}
q_{82} > M> n_2> k_2-k_1 \quad \text{and} \quad a(M):= \max\{a_i: 0\le i\le 82\} = a_{12}=134,
\end{eqnarray*}
from which it follows that $ m_3\le 100 $. We now return to \eqref{DAM1} and we put $ X:=1.12\times 10^{40} $ and $ C:=(20X)^5 $ and then apply the LLL algorithm in Lemma \ref{LLL} to $ m_3\in[1, 100] $. After a computer search in \textit{Mathematica}, we get
\begin{eqnarray*}
1.04\times 10^{-139}< \min_{1\le m_3\le 100}|\Gamma_4|< 24k_2\alpha^{-2\min\{n_i, m_j\}},
\end{eqnarray*}
then $ \min\{n_i, m_j\}\le 410 $. By continuing under the assumption that $ m_j:=\min\{n_i, m_j\} \le 426 $, we return to \eqref{DAM2} and put $ X:=1.12\times 10^{40} $, $ C:=(20X)^5 $ and $ M:=1.12\times 10^{38} $ for the case $ m_3< m_4 $ and the case $ m_3=m_4 $. After a computer search, we confirm that
\begin{eqnarray}
4.39\times 10^{-168} < \min_{\substack{1\le m_3\le 100 \\ m_3+1\le m_4 \le 426}} |\Gamma_4| < 24k_2\alpha^{-2n_3}.
\end{eqnarray}
This gives $ n_1\le 494 $ which holds in both cases. Hence, by a similar procedure given in the first cycle, we get that $ n_2 < 3\times 10^{36} $. 

We record what we have proved.
\begin{lemma}\label{firstredn1}
Let $ (k_i, n_i, m_i\} $ be a solution to the Diophantine equation $ x_{k_i}=L_{n_i}L_{m_i} $, with $ 0\le m_i\le  n_i $ for $ i=1,2 $ and $ 1\le k_1\le k_2 $, then
\begin{eqnarray*}
\max\{k_1, m_1\} \le  n_1 \le 494 \quad \text{and} \quad \max\{k_2, m_2\}\le  n_2 < 3\times 10^{36}.
\end{eqnarray*}
\end{lemma}

\subsection{The final reduction}
Returning back to \eqref{Pelling12} and \eqref{DM4} and using the fact that $ (x_1, y_1) $ is the smallest positive solution to the Pell equation \eqref{Pelleqn}, we obtain
\begin{eqnarray*}
x_k &=& \dfrac{1}{2}(\delta^k+\eta^k)  = \dfrac{1}{2}\left(\left(x_1+y_1\sqrt{d}\right)^{k}+\left(x_1-y_1\sqrt{d}\right)^k\right)\\
&=&\dfrac{1}{2}\left(\left(x_1+\sqrt{x_1^2\mp 1}\right)^{k}+\left(x_1-\sqrt{x_1^2\mp 1}\right)^k\right): = P^{\pm}_{k}(x_1).
\end{eqnarray*}
Thus, we return to the Diophantine equation $ x_{k_1}=L_{n_1}L_{m_1} $ and consider the equations 
\begin{eqnarray}\label{water1}
P^{+}_{k_1}(x_1)=L_{n_1}L_{m_1} \quad \text{and} \quad P^{-}_{k_1}(x_1)=L_{n_1}L_{m_1},
\end{eqnarray}
with $ k_1\in [1, 500] $, $ m_1\in [0, 500] $ and $n_1\in[m_1+1, 500] $.

Besides the trivial case $ k_1=1 $, with the help of a computer search in \textit{Mathematica} on the above equations in \eqref{water1}, we list the only nontrivial solutions in Table \ref{tab1} below. We also note that $$ 7+5\sqrt{2}= (1+\sqrt{2})^3, $$ so these solutions come from the same Pell equation with $ d=2 $.

\begin{table}[H]
\parbox{.47\linewidth}{
\centering
\begin{tabular}{|c c c c c|}
\hline
& & &$ Q^{+}_{k_1}(x_1) $&\\
\hline
$ k_1 $& $ x_1 $& $ y_1 $& $ d $ & $ \delta $\\
\hline
$ 2 $&$ 2 $& $ 1 $&$ 3 $ & $ 2+\sqrt{3} $\\
$ 2 $&$ 5 $& $ 2 $&$ 6 $ & $ 5+2\sqrt{6} $\\
$ 2 $&$ 10 $& $ 3 $&$ 11 $ & $ 10+3\sqrt{11} $\\
$ 2 $&$ 4 $& $ 1 $&$ 15 $ & $ 4+\sqrt{15} $\\
$ 2 $&$ 6 $& $ 1 $&$ 35 $ & $ 6+\sqrt{35} $\\
\hline
\end{tabular}
}
\hfill
\parbox{.47\linewidth}{
\centering
\begin{tabular}{|c c c c c|}
\hline
& & &$ Q^{-}_{k_1}(x_1) $&\\
\hline
$ k_1 $& $ x_1 $& $ y_1 $& $ d $ & $ \delta $\\
\hline
$ 2 $&$ 1 $& $ 1 $&$ 2 $ & $ 1+\sqrt{2} $\\
$ 2 $&$ 2 $& $ 1 $&$ 5 $ & $ 2+\sqrt{5} $\\
$ 2 $&$ 7 $& $ 5 $&$ 2 $ & $ 7+5\sqrt{2} $\\
$ 2 $&$ 4 $& $ 1 $&$ 17 $ & $ 4+\sqrt{17} $\\
$ 2 $&$ 26 $& $ 1 $&$ 677 $ & $ 26+\sqrt{677} $\\
$ 2 $&$ 179 $& $ 1 $&$ 32042 $ & $ 179+\sqrt{32042} $\\
\hline
\end{tabular}
}
\caption{Solutions to $ P^{\pm}_{k_1}(x_1)=L_{n_1}L_{m_1} $ }\label{tab1}
\end{table}
From the above tables, we set each $ \delta:=\delta_{t} $ for $ t=1, 2, \ldots 10 $.
We then work on the linear forms in logarithms $ \Gamma_1 $ and $ \Gamma_2 $, in order to reduce the bound on $ n_2 $ given in  Lemma \ref{firstredn1}. From the inequality \eqref{Pelling13}, for $ (k,n,m):=(k_2, n_2, m_2) $, we write
\begin{eqnarray}\label{mukazi111}
\left|k_2\dfrac{\log\delta_t}{\log\alpha}-(n_2+m_2)+\dfrac{\log 2}{\log (\alpha^{-1})}\right|<\left(\frac{12}{\log \alpha}\right)\alpha^{-2m_2}, 
\end{eqnarray}
for $t=1,2, \ldots 10$.

We put
\begin{eqnarray*}
\tau_{t}:=\dfrac{\log\delta_t}{\log\alpha}, \qquad \mu_t:=\dfrac{\log 2}{\log(\alpha^{-1})}\qquad \text{and} \quad (A_t, B_t):=\left(\frac{12}{\log \alpha}, \alpha\right).
\end{eqnarray*}
We note that $ \tau_t $ is transcendental by the Gelfond-Schneider's Theorem and thus, $ \tau_t $ is irrational. We can rewrite the above inequality, \eqref{mukazi111} as
\begin{eqnarray}\label{mukaziii11}
0<|k_2\tau_t-(n_2+m_2)+\mu_t|<A_tB_t^{-2m_2}, \quad \text{for} \quad t=1, 2, \ldots, 10.
\end{eqnarray}
We take $ M:= 3\times 10^{36} $ which is the upper bound on $ n_2 $ according to Lemma \ref{firstredn1} and apply Lemma \ref{Dujjella} to the inequality \eqref{mukaziii11}. As before, for each $ \tau_t $ with $ t=1, 2, \ldots, 10 $, we compute its continued fraction $ [a_0^{(t)}, a_1^{(t)}, a_2^{(t)}, \ldots ] $ and its convergents $ p_0^{(t)}/q_0^{(t)}, p_1^{(t)}/q_1^{(t)}, p_2^{(t)}/q_2^{(t)}, \ldots $. For each case, by means of a computer search in \textit{Mathematica}, we find and integer $s_{t}$ such that
\begin{eqnarray*}
q^{(t)}_{s_t}> 18\times 10^{36}=6M \qquad \text{ and } \qquad \varepsilon_t:=||\mu_t q^{(t)}||-M||\tau_t q^{(t)}|>0.
\end{eqnarray*}
We finally compute all the values of $ b_t:=\lfloor \log(A_t q^{(t)}_{s_t}/\epsilon_t)/\log B_t \rfloor /2$. The values of  $ b_t $ correspond to the upper bounds on $ \displaystyle{ m_2 } $, for each $ t=1, 2, \ldots, 10 $, according to Lemma \ref{Dujjella}.

Note that we have a problem at $ \delta_7 := 2+\sqrt{5} $. This is because $$ 2+\sqrt{5} = 2 \left(\frac{1+\sqrt{5}}{2}\right)^2 = 2\alpha^2. $$
So in this case we have $\Gamma_{1}:= (k_2-1)\log2 -(n_2+m_2-2k_2)\log\alpha$. Thus,
\begin{eqnarray*}
\left|\dfrac{\log2}{\log\alpha} - \dfrac{n_2+m_2-2k_2}{k_2-1}\right| < \dfrac{12}{(k_2-1)\alpha^{2m_2}\log\alpha}
\end{eqnarray*}
By a similar procedure given in Subsection \ref{RED} with $M:=3\times 10^{36}$, we get that $ q_{77}> M  $ and  $ a(M):=\max\{a_i: 0\le i \le 77\}=134   $. From this we can conclude that $ m_2\le 96 $.

The results of the computation for each $ t $ are recorded in Table \ref{tab2} below.
\begin{table}[H]
\begin{center}
\begin{tabular}{lllccc}
\hline 
$t$ & $\delta_t$ & $s_t$ & $q_{s_t} $ & $\varepsilon_t>$ & $b_t$ \\ 
\hline 
$1$ & $2+\sqrt{3}$ & $68$ & $ 2.07577\times 10^{37} $ & $ 0.319062 $ & $94$ \\ 
$2$ & $ 5+2\sqrt{6} $ & $ 91 $ & $ 8.19593\times 10^{37} $ & $ 0.087591 $ & $97$ \\  
$3$ & $ 10+3\sqrt{11} $& $ 67 $ & $ 2.25831\times 10^{38} $ & $ 0.316767 $ & $96$ \\  
$4$ & $4+\sqrt{15} $ & $ 70 $ & $  2.78896 \times 10^{37}$ & $ 0.329388 $ & $94$ \\  
$5$ & $6+\sqrt{35} $ & $ 74 $ & $ 1.75745\times 10^{38} $ & $ 0.409752 $ & $96$ \\ 
$6$ & $1+\sqrt{2}$ & $76$ & $ 2.02409\times 10^{37} $ & $ 0.263855 $ & $94$ \\ 
$7$ & $ 2+\sqrt{5} $ & $ - $ & $ - $ & $ - $ & $96$ \\  
$8$ & $ 4+\sqrt{17} $& $ 78 $ & $ 4.76137\times 10^{37} $ & $ 0.131771 $ & $96$ \\  
$9$ & $26+\sqrt{677} $ & $ 65 $ & $  3.17521 \times 10^{37}$ & $ 0.356148 $ & $94$ \\  
$10$ & $179+\sqrt{32042} $ & $ 77 $ & $ 3.45317\times 10^{37} $ & $ 0.384127 $ & $94$ \\ 
\hline 
\end{tabular} 
\end{center}
\caption{First reduction computation results}\label{tab2}
\end{table}

By replacing $ (k, n, m):=(k_2, n_2, m_2) $ in the inequality \eqref{DM4}, we can write
\begin{eqnarray}\label{water112}
\left|k_2\dfrac{\log\delta_t}{\log\alpha}-n_2+\dfrac{\log(2L_{m_2})}{\log(\alpha^{-1})}\right|< \left(\dfrac{12}{\log\alpha}\right)\alpha^{-2n_2}, 
\end{eqnarray}
for  $t=1,2,\ldots, 10$.

We now put
\begin{eqnarray*}
\tau_{t}:=\dfrac{\log\delta_t}{\log\alpha}, \quad \mu_{t, m_2}:=\dfrac{\log(2L_{m_2})}{\log(\alpha^{-1})}\quad \text{and} \quad (A_t, B_t):=\left(\frac{12}{\log \alpha}, \alpha\right).
\end{eqnarray*}
With the above notations, we can rewrite \eqref{water112} as
\begin{eqnarray}\label{water21}
0<|k_2\tau_t - n_2+\mu_{t, m_2}|<A_tB_t^{-2n_2}, \quad \text{ for} \quad t=1,2, \ldots 10.
\end{eqnarray}
We again apply Lemma \ref{Dujjella} to the above inequality \eqref{water21}, for 
\begin{eqnarray*}
t=1, 2, \ldots, 10, \quad m_2 =1, 2, \ldots, b_t, \quad \text{with}\quad M:=3\times 10^{36}.
\end{eqnarray*}
We take
\begin{eqnarray*}
\varepsilon_{t, m_2}:=||\mu_t q^{(t, m_2)}|| -M||\tau_t q^{(t, m_2)}||>0,
\end{eqnarray*}
and 
\begin{eqnarray*}
b_t=b_{t, m_2}:=\lfloor \log(A_t q^{(t, m_2)}_{s_t}/\epsilon_{t, m_2})/\log B_t \rfloor/2.
\end{eqnarray*}
The case $ \delta_7=2+\sqrt{5} $ is again treated individually by a similar procedure as in the previous step. With the help of \text{Mathematica}, we record the results of the computation in Table \ref{tab3} below.
\begin{table}[H]
\begin{center}
\begin{tabular}{c|cccccccccc}
\hline
$t$& $ 1 $&$ 2 $&$ 3 $&$ 4 $&$ 5 $& $ 6 $& $ 7 $& $ 8 $& $ 9 $&$ 10 $\\
$\varepsilon_{t,m_2}>$& $ 0.0145 $&$ 0.0002 $&$ 0.0006 $&$ 0.0034 $&$ 0.0106 $& $ 0.0005 $& $ - $& $ 0.0009 $& $ 0.0019 $&$ 0.0010  $\\
$b_{t, m_2}$& $ 97 $&$ 103 $&$ 102 $&$ 99 $&$ 99 $& $ 100 $ & $ 102 $& $ 100 $& $ 99 $& $ 100 $\\
\hline
\end{tabular}
\caption{Final reduction computation results}\label{tab3}
\end{center}
\end{table}
\begin{eqnarray*}
\text{Therefore, }~\max\{b_{t, m_2}: t=1, 2, \ldots, 10\quad \text{and} \quad m_2 = 1, 2, \ldots b_t\} \le 103 .
\end{eqnarray*}
Thus, by Lemma \ref{Dujjella}, we have that $ n_2\le 103 $, for all $ t=1,2, \ldots, 10 $. From the fact that $ \delta^{k}\le \alpha^{n+m+6} $, we can conclude that $ k_1< k_2\le 198  $. Collecting everything together, our problem is reduced to search for the solutions for \eqref{Problem} in the following ranges
\begin{eqnarray*}
1\le k_1<k_2\le 200, \quad 0\le m_1\le n_1 \le  200 \quad \text{and} \quad 0\le m_2\le n_2\le  200.
\end{eqnarray*}
After a computer search on the equation \eqref{Problem} on the above ranges, we obtained the following solutions, which are the only solutions for the  exceptional $ d $ cases  we have stated in Theorem \ref{Main}:

For the  $ +1 $ case:
\begin{eqnarray*}
(d=3)&& x_1 =2=L_1L_0,\quad  x_2=7=L_4L_1;\\
(d=15)&& x_1 =4=L_3L_1 =L_0L_0, \quad  x_5 =15124=L_{11}L_9; \\
(d=35)&& x_1=6 = L_{2}L_0, \quad x_3=846 = L_8L_6. 
\end{eqnarray*}

For the $ -1 $ case:
\begin{eqnarray*}
(d=2) && x_1=1=L_3L_3, \quad x_2=3=L_2L_1, \quad x_3=7=L_4L_1, \quad x_9=1393 = L_{11}L_4;\\
(d=5)&& x_1=2=L_1L_0, \quad x_2=9=L_{2}L_{2};\\
(d=17)&& x_1=4=L_3L_1=L_0L_0, \quad x_2=33=L_{5}L_{2}.
\end{eqnarray*}
This completes the proof of Theorem \ref{Main}. \qed

\section*{Acknowledgements}The author would like to thank the referee for the careful reading of the manuscript and the useful comments and suggestions that greatly improved on the quality of the presentation of this paper. He also thanks the referee, in particular for his/her contribution to the proofs of Lemma \ref{musu1} and Lemma \ref{musu2}. The author was supported by the Austrian Science Fund (FWF) grants: F5510-N26 -- Part of the special research program (SFB), ``Quasi-Monte Carlo Methods: Theory and Applications'' and W1230 --``Doctoral Program Discrete Mathematics''.

\medskip

\noindent MSC2010: 11B39, 11D45, 11D61, 11J86.

\end{document}